\documentclass[11pt]{article}
\usepackage{color}
\usepackage[latin1]{inputenc}
\usepackage[T1]{fontenc}
\usepackage[english]{babel}
\usepackage{amsfonts}
\usepackage{amssymb}
\usepackage{graphicx}
\usepackage{dsfont}
\usepackage{epsfig}
\usepackage{ae}
\usepackage{hyperref}
\usepackage{harvard}
\usepackage{amsmath}
\usepackage{setspace}
\usepackage{amsthm}
\usepackage{appendix}

\renewcommand{\(}{\left(}
\renewcommand{\)}{\right)}

\newcommand{\1}{\mathds{1}}

\renewcommand{\proof}{{\noindent \it Proof. }}

\newtheorem{theo}{Theorem}[section]
\newtheorem{pr}{Proposition}[section]
\newtheorem{lem}{Lemma}[section]
\newtheorem{co}{Corollary}[section]
\newtheorem{re}{Remark}[section]

\newcommand{\be}{\begin{eqnarray}}
\newcommand{\ee}{\end{eqnarray}}
\newcommand{\by}{\begin{eqnarray*}}
\newcommand{\ey}{\end{eqnarray*}}
\newcommand{\bt}{\begin{theo}}
\newcommand{\et}{\end{theo}}
\newcommand{\bl}{\begin{lem}}
\newcommand{\el}{\end{lem}}
\newtheorem*{lem1}{Lemma A.1}

\newtheorem*{lem2}{Lemma A.2}

\newtheorem*{lem3}{Lemma A.3}

\newcommand{\bc}{\begin{co}}
\newcommand{\ec}{\end{co}}
\newcommand{\eex}{\end{exa}\vspace{-3mm}}
\newcommand{\br}{\begin{re}}
\newcommand{\er}{\end{re}\vspace{-3mm}}
\renewcommand{\geq}{\geqslant}
\renewcommand{\leq}{\leqslant}

\citationstyle{agsm}

\begin{document}

\title{A new proof of an Engelbert-Schmidt type zero-one law for time-homogeneous diffusions\footnote{The author thanks two anonymous referees for helpful suggestions in improving the presentation of the paper.  The author is also grateful to Carole Bernard and Don McLeish for helpful discussions. The usual disclaimer applies.}}
\author{   \textsc{Zhenyu Cui} \footnote{Corresponding Author. Zhenyu
Cui is with the department of Mathematics at Brooklyn College of the City University of New York, Ingersoll Hall, 2900 Bedford Ave, Brooklyn, NY11210, United States.   Tel.: +1718-951-5600, ext. 6892
              Fax: +1718-951-4674.
              Email: \texttt{zhenyucui@brooklyn.cuny.edu} } \ }

\date{Draft: \today }
\maketitle
\begin{abstract}
In this paper we give a new proof to an Engelbert-Schmidt type zero-one law for time-homogeneous diffusions, which provides deterministic criteria for the
convergence of integral
functional of diffusions. Our proof is based on a slightly stronger assumption than that in Mijatovi\'c and Urusov \citeyear{MU12WP}, and utilizes stochastic time change and Feller's test of explosions. It does not rely on advanced methods such as the first Ray-Knight theorem, Wiliam's theorem, Shepp's dichotomy result for Gaussian processes or Jeulin's lemma as in the previous literature(see Mijatovi\'c and Urusov \citeyear{MU12WP} for a pointer to the literature). The new proof has an intuitive interpretation as we link the integral functional to the explosion time of an associated diffusion process.
\end{abstract}


\newpage

\section{Introduction}

The Engelbert-Schmidt zero-one law was initially proved in the standard Brownian motion case (see Engelbert and Schmidt \citeyear{ES81} or Proposition $3.6.27$, p$216$ of  Karatzas and Shreve \citeyear{KS91}). Engelbert and Tittel \citeyear{ET02} obtain a generalized Engelbert-Schmidt type zero-one law for the integral functional $\int_0^t f(X_s)ds$, where $f$ is a non-negative Borel measurable function and $X$ is a strong Markov continuous local martingale. In an expository paper, Mijatovi\'c and Urusov \citeyear{MU12WP} consider the case of a one-dimensional time-homogeneous diffusion and their Theorem $2.11$ gives the corresponding Engelbert-Schmidt type zero-one law for time-homogeneous diffusions. Their proof relies on their Lemma $4.1$, which they provide two proofs without using the Jeulin's lemma(see Lemma $3.1$ of Engelbert and Tittel \citeyear{ET02}). The first proof is based on William's theorem (Ch.VII, Corollary $4.6$, p$317$, Revuz and Yor \citeyear{RY99}). The second proof is based on the first Ray-Knight theorem (Ch.XI, Theorem $2.2$, p$455$, Revuz and Yor \citeyear{RY99}). In Khoshnevisan, Salminen and Yor \citeyear{KSY06}, a similar question on the convergence of integral functional of diffusions is treated under different assumptions with answers given in different terms (see the discussion on p2 of Mijatovi\'c and Urusov \citeyear{MU12WP}). Here our discussion is based on the setting in Mijatovi\'c and Urusov \citeyear{MU12WP}.

The contribution of this paper is two-fold. First, under a slightly stronger assumption that $f$ is positive, we complement the study of the Engelbert-Schmidt type zero-one law in Mijatovi\'c and Urusov \citeyear{MU12WP} with a new simple proof that circumvents advanced tools such as the William's theorem, and the first Ray-Knight theorem, which are employed by them to prove their Lemma 4.1. As discussed on p10 in their paper, in previous literature, this lemma has been proven using the first Ray-Knight theorem with either  Shepp's dichotomy result(Shepp \citeyear{S66}) for Gaussian processes(see Engelbert and Schmidt \citeyear{ES87}) or the Jeulin's lemma(see Engelbert and Tittel \citeyear{ET02}).
Second, through stochastic time change, we establish a link between the integral functional of diffusions and the explosion time of another associated time-homogeneous diffusion.  Our proof has an intuitive interpretation: ``\textit{the convergence/divergence of integral functional of diffusions is equivalent to the explosion/non-explosion of the associated diffusion}". Mijatovi\'c and Urusov \citeyear{MU12WP} give a proof(see their Proposition $2.12$) of the Feller's test of explosions as an application of their results, and Theorem \ref{tczo1} here provides a converse to their result under a slightly stronger assumption.

The paper is organized as follows. Section \ref{s2} describes the probabilistic setting, states results on stochastic time change and gives the new proof. Section \ref{s3} concludes the paper and provides future research directions.

\section{Main Result\label{s2}}

Given a complete filtered probability space $(\Omega, \mathcal{F}, (\mathcal{F}_t)_{t\in[0,\infty)}, P)$ with state space $J=(\ell,r), -\infty\leq \ell<r\leq \infty$,  and $\mathcal{F}$ is right continuous(i.e. $\mathcal{F}_t=\mathcal{F}_{t+}$ for $t\in[0,\infty)$). Assume that the $J$-valued diffusion $Y=(Y_t)_{\left\{t\in[0,\infty)\right\}}$ satisfies the stochastic differential equation(SDE)
\begin{align}
dY_{t} &=\mu (Y_t) dt +\sigma(Y_t)  dW_t,\quad Y_0=x_0,\label{y}
\end{align}
where $W$ is a $\mathcal{F}_t$-Brownian motion and $\mu, \sigma: J\rightarrow \mathbb{R}$ are Borel measurable functions satisfying the Engelbert-Schmidt conditions
\begin{align}
\forall x\in J, \text{ } \sigma(x)\neq 0,\quad  \frac{1}{\sigma^2(\cdot)}, \quad \frac{\mu(\cdot)}{\sigma^2(\cdot)}\in L_{loc}^1 (J),\label{cond1}
\end{align}
where $L_{loc}^{1} (J)$ denotes the class of locally integrable functions, i.e. the functions $J\rightarrow \mathbb{R}$ are integrable on compact subsets of $J$. This condition \eqref{cond1} guarantees that the SDE \eqref{y} has a unique in law weak solution that possibly exits its state space $J$(see Theorem $5.15$, p$341$, Karatzas and Shreve \citeyear{KS91}).

 Denote the possible explosion time of $Y$ from its state space by $\zeta$, i.e. $\zeta=\inf\{u>0, Y_u \not\in J\} $,
  which means that $P$-a.s. on $\{\zeta=\infty\}$ the trajectories of $Y$ do not exit $J$, and $P$-a.s. on $\{\zeta<\infty\}$,
  we have $\lim\limits_{t\rightarrow \zeta}Y_t =r$ or $\lim\limits_{t\rightarrow \zeta}Y_t =\ell$. $Y$ is defined such that it stays at its
  exit point, which means that $l$ and $r$ are absorbing boundaries. The following terminology is used: $Y$ exits the state space $J$ at $r$
  means $P\(\zeta<\infty, \lim\limits_{t\rightarrow \zeta}Y_t=r\)>0$.

Let $b$ be a Borel measurable function such that $\forall x\in J, b(x)\neq 0$, and assume the following local integrability condition
\begin{align}
   \frac{b^2(\cdot)}{\sigma^2(\cdot)} \in L_{loc}^1 (J).\label{cond2}
\end{align}
Define the function  $\varphi_t:=\int_0^t b^2(Y_u)du, \text{ for }t\in [0,\zeta]$.
\begin{lem}\label{ct}
Assume \eqref{cond2} and also that  $\forall x\in J, b(x)\neq 0$, then $\varphi_t$ is a strictly increasing  function for $t\in [0,\zeta]$ and it is absolutely continuous for a.a. $\omega$'s on compact intervals of $[0,\zeta)$. Furthermore, $\varphi_t <\infty$ for $t\in [0,\zeta)$
\end{lem}
\proof
Since $\forall x\in J, b(x)\neq 0$, $\varphi_t$ is a strictly
 increasing function for $t\in [0,\zeta]$. Note that $\varphi_t, t\in [0, \zeta)$ is represented as a time integral, and the continuity follows.  It is a standard result that the condition \eqref{cond2} implies that $\varphi_t <\infty$ for $t\in [0,\zeta)$(see the proof after (8) on p4 and p5 in Mijatovi\'c and Urusov \citeyear{MU12PTRF}).
\hfill$\Box$\\

 The following result is about stochastic time change.
\begin{pr}\label{tcm}(Theorem $3.2.1$ of Cui \citeyear{C13})

Recall that $\varphi_t=\int_0^t b^2(Y_u)du, t\in [0,\zeta]$. Assume that the conditions \eqref{cond1}, \eqref{cond2} are satisfied, and $\forall x\in J, b(x)\neq 0$.

(i) Under $(\Omega, \mathcal{F}, (\mathcal{F}_t)_{t\in[0,\infty)}, P)$, define
 \begin{align}
 T_t&:=
 \begin{cases}
 \inf\{u\geq 0 : \varphi_{u\wedge \zeta} >t\}, &\mbox{ on } \left\{0\leq t <\varphi_{\zeta}\right\},\\
 \infty, &\mbox{ on } \left\{\varphi_{\zeta}\leq t<\infty\right\}.\label{ttc}
 \end{cases}
 \end{align}
 Define a new filtration $\mathcal{G}_{t}:=\mathcal{F}_{T_t}, t\in[0,\infty)$, and a new process $X_t:=Y_{T_t}$, on $\left\{0\leq t <\varphi_{\zeta}\right\}$. Then $X_t$ is $\mathcal{G}_t$-adapted and we have the stochastic representation
\begin{align}
Y_{t} &=X_{\int_{0}^{t} b^2(Y_{s})ds}=X_{\varphi_{t}}, \text{ $P$-a.s. }\quad \text{ on } \left\{0\leq t <\zeta \right\},\label{tc1fs}
\end{align}
and the process $X$ is a time-homogeneous diffusion, which solves the following SDE under $P$
\begin{align}
dX_{t} &=\frac{\mu(X_t)}{b^2(X_t)}\1_{t\in [0,\varphi_{\zeta})}dt +\frac{\sigma(X_t)}{b(X_t)}\1_{t\in [0,\varphi_{\zeta})}dB_t, \quad X_0=x_0,\label{xsde1}
\end{align}
where $B_t$ is the $\mathcal{G}_t$-adapted Dambis-Dubins-Schwartz Brownian motion  under $P$ defined in the proof. 

(ii) Define $\zeta^X :=\inf\left\{u>0: X_u \not \in J \right\}$, then $\zeta^X =\varphi_{\zeta}=\int_0^{\zeta} b^2 (Y_s )ds$, $P$-a.s., and we can rewrite the SDE \eqref{xsde1} under $P$ as
\begin{align}
dX_{t} &=\frac{\mu(X_t)}{b^2(X_t)}\1_{t\in [0,\zeta^X)}dt +\frac{\sigma(X_t)}{b(X_t)}\1_{t\in [0,\zeta^X)}dB_t, \quad X_0=x_0.\label{xsdec}
\end{align}

(iii) The event $\left\{\limsup\limits_{t\rightarrow \zeta} Y_t =r\right\}$ is equal to $\left\{\limsup\limits_{t\rightarrow \zeta^{X}} X_t =r\right\}$.
Similarly for the case of the left boundary $\ell$, the case of $\liminf$, and $\lim$. 
\end{pr}

\begin{proof}
From Lemma \ref{ct}, $\varphi_s$ is a strictly increasing function on $[0,\zeta]$. From Problem $3.4.5$ (ii), p174 of Karatzas and Shreve \citeyear{KS91},  $\varphi_{T_t \wedge \zeta}=t\wedge \varphi_{\zeta}$, $P$-a.s. for $0\leq t<\infty$. On $\left\{0\leq t<\varphi_{\zeta}\right\}$, take $u=\zeta$ in \eqref{ttc},  then $\varphi_{\zeta\wedge \zeta}=\varphi_{\zeta}>t$ holds $P$-a.s. according to the assumption. Then  $T_t\leq \zeta$, $P$-a.s. because of the definition \eqref{ttc}, $T_t:=\inf\{u\geq 0 : \varphi_{u\wedge \zeta} >t\}$. Thus  $\varphi_{T_t}=t$, $P$-a.s. on $\left\{0\leq t<\varphi_{\zeta}\right\}$.

Choose $t=\varphi_s$ on $\left\{0\leq s<\zeta\right\}$, then $0\leq t <\varphi_{\zeta}$, $P$-a.s. After substituting this $t$ into the definition of the process $X$, we have $X_{\varphi_s}=X_t:=Y_{T_t}=Y_{T_{\varphi_s}}=Y_s$, $P$-a.s.. For the last equality, recall the definition and  $T_{\varphi_s}=\inf\{u\geq 0:\varphi_{u\wedge\zeta}>\varphi_s \}=\inf\{u\geq 0:u\wedge\zeta>s \}=s$, $P$-a.s., on $\left\{0\leq s<\zeta\right\}$. Then we have proved the representation $Y_s =X_{\varphi_s}$, $P$-a.s. on $\left\{0\leq s<\zeta\right\}$, and the next goal is to determine the coefficients of the SDE satisfied by $X$ under $P$.

For $X$ satisfying the relation \eqref{tc1fs}, we aim to show that $X$ satisfies the following SDE under $P$ \begin{align}
dX_{t} &=\frac{\mu(X_t)}{b^2(X_t)}\1_{t\in [0,\varphi_\zeta)}dt +\frac{\sigma(X_t)}{b(X_t)}\1_{t\in [0,\varphi_\zeta)}dB_t,\quad  X_0=Y_0=x_0. \label{x1c}
\end{align}
where $B$ is the Dambis-Dubins-Schwartz Brownian motion adapted to $\mathcal{G}_t$ constructed as follows:

Note that $M_{t\wedge \zeta}=\int_0^{t\wedge \zeta} b(Y_u) dW_u, t\in[0,\infty)$ is a continuous local martingale with quadratic variation
$\varphi_{t\wedge \zeta}=\int_0^{t\wedge \zeta} b^2(Y_u) du, t\in[0,\infty)$.
Then $\lim\limits_{t\rightarrow \infty}\varphi_{t\wedge \zeta}=\varphi_{\zeta}$, $P$-a.s. due to the  continuity of $\varphi_s$(see Lemma \ref{ct}).

From the Dambis-Dubins-Schwartz theorem (Ch.V, Theorem $1.6$ and Theorem $1.7$ of Revuz and Yor \citeyear{RY99}),  there exists
an enlargement $(\overline{\Omega},\overline{\mathcal{G}_t})$ of $(\Omega, \mathcal{G}_{t})$ and a standard Brownian motion $\overline{\beta}$ on $\overline{\Omega}$ independent of
$M$ with $\overline{\beta}_0=0$, such that the process
\begin{align}
B_{t} &:=
\begin{cases}
\int_0^{T_t} b(Y_u) dW_u, &\mbox{ on } \left\{t < \varphi_{\zeta}\right\}, \\
\int_0^{\zeta} b(Y_u) dW_u+\widetilde{\beta}_{t-\varphi_{\zeta}}, &\mbox{ on } \left\{t \geq \varphi_{\zeta}\right\}.\label{conbc}
\end{cases}
\end{align}
is a standard linear Brownian motion. Our construction of $T_t, t\in[0,\infty)$ agrees with that in Problem $3.4.5$, p174 of Karatzas and Shreve \citeyear{KS91}. From Problem $3.4.5$ (ii) and the construction \eqref{conbc},  $B_{\varphi_s}=M_s$, $P$-a.s. on $\left\{0\leq s<\zeta\right\}$, and on $\left\{s=\zeta\right\}$, $B_{\varphi_{\zeta}}:=\int_0^{\zeta} b(Y_u) dW_u=:M_{\zeta}$, $P$-a.s.. Thus $B_{\varphi_t}=M_t$, $P$-a.s. on $\left\{0\leq t\leq \zeta\right\}$.

For the convenience of exposition, denote $\mu_1(.)=\mu(.)/b^2(.)$, and $\sigma_1(.)=\sigma(.)/b(.)$.
Integrate the SDE in \eqref{y} under $P$ from $0$ to $t\wedge \zeta$
 \begin{align}
 Y_{t\wedge \zeta} -Y_0&=\int_0^{t\wedge \zeta} \mu(Y_u) du +\int_0^{t\wedge \zeta} \sigma(Y_u)dW_u\notag\\
 &= \int_0^{t\wedge \zeta} \mu_1(Y_u)b^2 (Y_u)du +\int_0^{t\wedge \zeta} \sigma_1(Y_u)b(Y_u)dW_u.\label{intermet2}
 \end{align}

Apply the change of variables formula similar to Problem $3.4.5$ (vi), p174 of Karatzas and Shreve \citeyear{KS91}, and note the relation \eqref{tc1fs}
\begin{align}
\int_0^{t\wedge \zeta} \mu_1(Y_{u})b^2(Y_{u})du&= \int_0^{t\wedge \zeta} \mu_1(X_{\varphi_{u}})d\varphi_{u}=\int_0^{\varphi_{t\wedge \zeta}} \mu_1(X_u)du, \label{pu1}
\end{align}
and similarly
\begin{align}
\int_0^{t\wedge \zeta} \sigma_1(Y_u)b(Y_{u})dW_u&=\int_0^{t\wedge \zeta} \sigma_1 (X_{\varphi_u}) dB_{\varphi_u}
= \int_0^{\varphi_{t\wedge \zeta}} \sigma_1 (X_u)dB_u, \label{pu2}
\end{align}
where the first equality in \eqref{pu2} is due to the relationship $B_{\varphi_u}=M_u=\int_0^u b(V_s)dW_s$, $P$-a.s., on $\left\{0\leq u \leq t\wedge\zeta\right\}$, which we have established above.
Also notice the representation $Y_{t\wedge \zeta}=X_{\varphi_{t\wedge \zeta}}$, $P$-a.s., and $Y_0=X_0$, then
\begin{align}
X_{\varphi_{t\wedge \zeta}} -X_0 &=\int_0^{\varphi_{t\wedge \zeta}} \mu_1(X_u)du +\int_0^{\varphi_{t\wedge \zeta}} \sigma_1 (X_u)dB_u. \label{intermet}
\end{align}
Then on $\left\{0\leq s\leq \varphi_{t\wedge \zeta}\right\}$
\begin{align}
X_s -X_0 &=\int_0^s \mu_1(X_u)du +\int_0^s \sigma_1 (X_u)dB_u.\label{a3}
\end{align}
Note that for $0\leq t<\infty$, we have $s\in[0, \varphi_{\zeta}]$, $P$-a.s. From \eqref{a3}, and recall the definition of $\mu_1(.)$ and $\sigma_1(.)$, we have the following SDE for $X$ under $P$:
\begin{align}
dX_{s} &=\frac{\mu(X_s)}{b^2(X_s)}\1_{s\in [0,\varphi_\zeta)}ds +\frac{\sigma(X_s)}{b(X_s)}\1_{s\in [0,\varphi_\zeta)}dB_s,\quad  X_0=Y_0=x_0. \notag
\end{align}
This completes the proof of statement (i).

Statement (ii) is a direct consequence of the stochastic representation $Y_{t\wedge \zeta}=X_{\varphi_{t\wedge \zeta}}$, $P$-a.s. in statement (i), because from Lemma \ref{ct}, $\varphi_t$ is a strictly increasing function in $t$.

For statement (iii), denote $f(t)=Y_t$ on $\left\{0\leq t<\zeta\right\}$ and $g(t)=X_t$ on $\left\{0\leq t<\zeta^X\right\}$. From statement (i),  $g(\varphi_t)=X_{\varphi_t}=Y_t=f(t)$, $P$-a.s. on $\left\{0\leq t<\zeta\right\}$. They are two real-valued functions linked by a strictly increasing and continuous function $\varphi_t$. From statement (ii),  $\varphi_{\zeta}=\zeta^{X}$, $P$-a.s.
 This means that
 $\limsup\limits_{t\rightarrow \zeta} Y_t =\limsup\limits_{t\rightarrow \zeta} f(t)=\limsup\limits_{t\rightarrow \zeta}g(\varphi_t)=\limsup\limits_{t\rightarrow \zeta^{X}}g(t)=\limsup\limits_{t\rightarrow \zeta^{X}}X_t$, $P$-a.s., and the equality of the two events holds. Similarly for the cases of $\liminf$ and $\lim$.   This completes the proof.
\hfill$\Box$\\
\end{proof}


Denote $\zeta_{\ell} (resp.  \zeta_r)$ as the possible explosion time of the diffusion $Y$ through the boundary $\ell  (resp.  r)$. Correspondingly, denote $\zeta^{X}_{\ell} (resp. \zeta^{X}_r)$ as the possible explosion time of the diffusion $X$ through the boundary $\ell  (resp. r)$.
Recall the definition $\zeta$, and we have $\zeta =\min (\zeta_{\ell},\zeta_r),  \zeta^{X}=\min( \zeta^{X}_{\ell},  \zeta^{X}_r)$.
From  Proposition \ref{tcm} (ii)
\begin{align}
\zeta^X_{\ell} &=\int_0^{\zeta_{\ell}} b^2 (Y_s )ds, \quad \zeta^X_r =\int_0^{\zeta_r} b^2 (Y_s )ds, \quad \zeta^X =\int_0^{\zeta} b^2(Y_s)ds, \text{ $P$-a.s.}\label{x2}
\end{align}

Fix an arbitrary constant $c\in J$ and introduce the scale function $s(.)$ of the SDE \eqref{y} under $P$
 \begin{align}
 s(x)&:=\int_c^x \exp\left\{-\int_{c}^{y} \frac{2\mu}{\sigma^2}(u)du  \right\}dy, \quad x\in \bar{J}.\label{scale}
 \end{align}
The scale function of the diffusion X with SDE \eqref{xsdec} is also $s(.)$.
With an arbitrary constant $c\in J$, for $x\in\bar{J}$, introduce the following test functions respectively for $Y$ and $X$
\begin{align}
v(x) \equiv \int_c^x (s(x)-s(y)) \frac{2}{s^{\prime}(y) \sigma^2(y)}dy, \quad
v_{X} (x) \equiv \int_c^x (s(x)-s(y)) \frac{2 b^2(y)}{s^{\prime}(y) \sigma^2(y)}dy.\notag\\
\label{v4}
\end{align}
Recall the classical Feller's test of explosions for diffusions using our notation.
\begin{lem}\label{f1}(Theorem $5.29$ of Karatzas and Shreve \citeyear{KS91})


For $Y$, $P(\zeta =\infty)=1$ if and only if $v(\ell)=v(r)=\infty$; Otherwise $P(\zeta <\infty)>0$.

For $X$, $P(\zeta^X =\infty)=1$ if and only if $v_X(\ell)=v_X(r)=\infty$; Otherwise $P(\zeta^X <\infty)>0$.

\end{lem}
The following result provides the precise conditions when $P(\zeta^X <\infty)=1$ holds.
\begin{lem}\label{f2} (Proposition $5.32$ of Karatzas and Shreve \citeyear{KS91})

We have that $P(\zeta^X <\infty)=1$  if and only if  at least one of the following statements are satisfied:

(a) $v_X(r)<\infty$ and $v_X({\ell})<\infty$;

(b) $v_X(r)<\infty$ and $s({\ell})=-\infty$;

(c) $v_X({\ell})<\infty$ and $s(r)=\infty$.

\end{lem}


From Feller's test of explosions, we have that the process $X$ under $P$  may exit its state space $J$ at the boundary point $r$, i.e. $P(\zeta^X <\infty, \lim\limits_{t\rightarrow \zeta^X} X_t =r )>0$, if and only if
\begin{align}
v_X (r)&<\infty.     \label{b}
\end{align}

Similarly for the case of the endpoint $\ell$.

We have the following exhaustive classifications of the possible explosion events of $Y$
\begin{align}
A&=\left\{\zeta=\infty, \limsup\limits_{t\rightarrow\infty} Y_t =r, \quad \liminf_{t\rightarrow \infty} Y_t =\ell\right\},\notag\\
B_r&=\left\{\zeta=\infty, \lim\limits_{t\rightarrow \infty}Y_t =r\right\}, \quad C_r =\left\{\zeta<\infty, \lim\limits_{t\rightarrow \zeta}Y_t =r\right\},\notag\\
B_{\ell} &=\left\{\zeta=\infty, \lim\limits_{t\rightarrow \infty} Y_t =\ell\right\}, \quad C_{\ell} =\left\{\zeta<\infty, \lim\limits_{t\rightarrow \zeta} Y_t =\ell\right\}.\notag
\end{align}
Similarly for $X$
\begin{align}
A^X&=\left\{\zeta^{X}=\infty, \limsup\limits_{t\rightarrow\infty} X_t =r, \quad \liminf\limits_{t\rightarrow \infty} X_t =\ell\right\},\notag\\
B_r^X&=\left\{\zeta^{X}=\infty, \lim\limits_{t\rightarrow \infty}X_t =r\right\}, \quad
C_r^X =\left\{\zeta^{X}<\infty, \lim\limits_{t\rightarrow \zeta^{X}} X_t =r\right\},\notag\\
B_{\ell}^X &=\left\{\zeta^{X}=\infty, \lim\limits_{t\rightarrow \infty} X_t =\ell\right\}, \quad
C_{\ell}^X =\left\{\zeta^{X}<\infty, \lim\limits_{t\rightarrow \zeta^{X}} X_t =\ell\right\}.\label{secondtc}
\end{align}
We first recall a lemma using our notation.
\begin{lem}\label{auxtc} (Proposition $2.3$, $2.4$ and $2.5$ on p4 of Mijatovi\'c and Urusov \citeyear{MU12WP} )

\vspace{0.1cm}

(1) Either $P(A^X)=1$ or $P(B_r^X \cup C_r^X \cup B_{\ell}^X \cup C_{\ell}^X)=1$.\\

(2) (i) $P(B_r^X \cup C_r^X)=0$ holds if and only if $s(r)=\infty$.\\

$\quad$ (ii) $P(B_{\ell}^X \cup C_{\ell}^X)=0$ holds if and only if $s(\ell)=-\infty$.\\

(3) Assume that $s(r)<\infty$. Then either $P(B_r^X)>0, P(C_r^X)=0$ or $P(B_r^X)=0, P(C_r^X)>0$. Similarly for the case of the endpoint $\ell$.
\end{lem}
\begin{proof}
For details of the proof,  refer to Mijatovi\'c and Urusov \citeyear{MU12WP} and Engelbert and Schmidt \citeyear{ES91}. \hfill$\Box$\\
\end{proof}

 Now we are ready to give the new proof to an  Engelbert-Schmidt type zero-one law for $Y$.
\begin{theo}\label{tczo1}(Engelbert-Schmidt type zero-one law for time-homogeneous diffusions, Theorem $2.12$ of Mijatovi\'c and Urusov \citeyear{MU12WP} with a stronger assumption\footnote{Theorem $2.12$ of Mijatovi\'c and Urusov \citeyear{MU12WP} assumes $f:J\rightarrow [0,\infty]$, which is weaker than the current assumption of $f$.})

Assume that the function $f:J\rightarrow (0,\infty]$ is a positive Borel measurable function and satisfies $f/\sigma^2 \in L_{loc}^1(J)$. Let $s(r)<\infty$.

(i)If $\frac{(s(r)-s)f}{s^{\prime}\sigma^2}\in L^1_{loc}(r-)$, then $\int_0^{\zeta} f(Y_u)du<\infty$, $P$-a.s. on $\left\{\lim\limits_{t\rightarrow \zeta} Y_t =r\right\}$.

(ii)If $\frac{(s(r)-s)f}{s^{\prime}\sigma^2}\not\in L^1_{loc}(r-)$, then $\int_0^{\zeta} f(Y_u)du=\infty$, $P$-a.s. on $\left\{\lim\limits_{t\rightarrow \zeta} Y_t =r\right\}$.

The analogous results on the set $\left\{\lim\limits_{t\rightarrow \zeta} Y_t =\ell\right\}$ can be similarly stated.
\end{theo}

\begin{proof}
To be consistent with our notation, let $b(x):=\sqrt
{f(x)}>0$ for $x\in J$. Then the assumptions of Lemma \ref{ct} and
Proposition \ref{tcm} are satisfied.

Denote $G=\left\{  \lim\limits_{t\rightarrow\zeta}
{Y}_{t}=r\right\}$, and from Proposition \ref{tcm} (iii)
\begin{align}
G&=\left\{  \lim\limits_{t\rightarrow\zeta}Y_{t}=r\right\}
=\left\{  \lim\limits_{t\rightarrow\zeta^{X}}X
_{t}=r\right\}  =B_{r}^{X}\cup C_{r}^{X}.\notag
\end{align}
 The result is trivial in the case $P(G)=0$,  so we assume
$P(G)>0.$  Since the events $B_{r}^{X},C_{r}^{X}$  are
disjoint
\begin{equation}
P(G)=P(B_{r}^{X})+P(C_{r}^{X}).\label{PG}%
\end{equation}
From Lemma \ref{auxtc}, ${s}(r)<\infty$ implies that either $ P(B_r^X)>0,  P(C_r^X)=0$ or $ P(B_r^X)=0,  P(C_r^X)>0$ holds.

For statement (i$^{\prime}$), $\frac
{( {s}(r)- {s})b^2}{ {s}^{\prime}\sigma^{2}}\in L_{loc}^{1}(r-)$, combined with $ {s}(r)<\infty$, implies $v_{X}(r)<\infty$.  From equation \eqref{b},  this is equivalent to $P\(\zeta^{X}<\infty, \lim\limits_{t\rightarrow \zeta^{X}} X_t =r \)>0$, and from \eqref{secondtc}, it means  $ P(C_r^X)>0$. Thus $ P(B_r^X)=0,  P(C_r^X)>0$ holds. This together with \eqref{PG} implies
\begin{align}
P(G)&= {P}(C_{r}^{X})=P\(\zeta^{X}<\infty, \lim\limits_{t\rightarrow \zeta} Y_t =r\)\notag\\
&=P\(\int_{0}^{\zeta}b^{2}(Y_{u})du<\infty, \lim\limits_{t\rightarrow \zeta} Y_t =r\)\notag\\
&=P\(\int_{0}^{\zeta}f(Y_{u})du<\infty, \lim\limits_{t\rightarrow \zeta} Y_t =r\),\notag
 \end{align}
where the third equality follows from Proposition \ref{tcm} (ii).

For statement (ii$^{\prime}$), $\frac{( {s}(r)- {s})b^2}{ {s}^{\prime}\sigma^{2}}\not \in L_{loc}^{1}(r-)$, combined with $ {s}(r)<\infty$, implies $ v_{X}(r)=\infty$. From equation \eqref{b},  this is equivalent to $P\(\zeta^{X}<\infty, \lim\limits_{t\rightarrow \zeta^{X}} X_t =r \)=0$, and from \eqref{secondtc}, it means  $ P(C_r^X)=0$. Thus $ P(B_r^X)>0,  P(C_r^X)=0$ holds. By a
similar argument to that above
\begin{align}
P(G)&=P(B_{r}^{X})=P\(\zeta^{X}=\infty, \lim\limits_{t\rightarrow \zeta} Y_t =r\)\notag\\
&=P\(\int
_{0}^{\zeta}b^{2}(Y_{u})du=\infty, \lim\limits_{t\rightarrow \zeta} Y_t =r\)\notag\\
&=P\(\int_{0}^{\zeta}f(Y_{u})du=\infty, \lim\limits_{t\rightarrow \zeta} Y_t =r\).\notag
\end{align}
The analogous results on the set $\left\{\lim\limits_{t\rightarrow
\zeta}Y_{t}=\ell\right\}  $ can be similarly proved by switching
the roles of $r$ and $\ell$ in the above argument. This completes the proof. \hfill$\Box$\\

\end{proof}

\section{Conclusion and future research\label{s3}}

In this paper, through stochastic time change, we have established a link between the integral functional of a diffusion and the explosion time of an associated time-homogeneous diffusion. A new proof to an Engelbert-Schmidt type zero-one law for diffusions is presented under a slightly stronger assumption than Mijatovi\'c and Urusov \citeyear{MU12WP}.
Recently Karatzas and Ruf \citeyear{KR13} give a detailed study of the distribution of this explosion time in a one-dimensional time-homogeneous diffusion setting, and present concrete  examples where the distribution function of the explosion time can be explicitly determined. Future research direction will be to study and search for explicit characterizations of the distribution functions of integral functional of diffusions.



\bibliography{zeroone}

\end{document}